\documentclass[12pt,a4paper]{article}
\usepackage{etoolbox} 
\usepackage{amsfonts}
\usepackage[margin=1in]{geometry}
\usepackage{times}
\usepackage{amsthm,amssymb}
\usepackage{amsmath,xypic}
\usepackage{multicol}
\usepackage{hyphenat}
\usepackage[usenames,dvipsnames]{color}
\usepackage{rotating}
\usepackage{longtable}
\usepackage{caption}
\usepackage{tcolorbox}

\usepackage{epigraph}
\usepackage{enumitem}
\setlist[enumerate]{listparindent=0.5in}
\usepackage{mathrsfs}
\DeclareMathAlphabet{\mathscrbf}{OMS}{mdugm}{b}{n} 
\newcommand{\be}{\begin{equation}}
\newcommand{\ee}{\end{equation}}
\newcommand{\bes}{\begin{equation*}}
\newcommand{\ees}{\end{equation*}}
\newcommand{\bea}{\begin{eqnarray}}
\newcommand{\eea}{\end{eqnarray}}
\newcommand{\beas}{\begin{eqnarray}}
\newcommand{\eeas}{\end{eqnarray}}
\newcommand{\ben}{\begin{note}}
\newcommand{\een}{\end{note}}
\newcommand{\bexl}{\vskip0.1em\noindent\hrulefill\vskip1em\begin{ExerciseList}}
\newcommand{\eexl}{\end{ExerciseList}\hrulefill}

\newcommand{\bthm}{\begin{theorem}}
\newcommand{\ethm}{\end{theorem}}
\newcommand{\bpro}{\begin{prop}}
\newcommand{\epro}{\end{prop}}
\newcommand{\bcor}{\begin{corollary}}
\newcommand{\ecor}{\end{corollary}}
\newcommand{\bcon}{\begin{conjecture}}
\newcommand{\econ}{\end{conjecture}}
\newcommand{\bp}{\begin{proof}}
\newcommand{\ep}{\end{proof}}
\newcommand{\blem}{\begin{lemma}}
\newcommand{\elem}{\end{lemma}}
\newcommand{\bn}{\begin{note}}
\newcommand{\en}{\end{note}}
\newcommand{\benum}{\begin{enumerate}}
\newcommand{\eenum}{\end{enumerate}}
\newcommand{\bed}{\begin{defn}}
\newcommand{\eed}{\end{defn}}
\newcommand{\brem}{\begin{remark}}
\newcommand{\erem}{\end{remark}}

\newcommand{\btik}{\begin{tikzpicture}\begin{axis}[scale=0.5,axis y line=center, axis x line=middle]}
\newcommand{\etik}{\end{axis}\end{tikzpicture}}

\let\into=\hookrightarrow
\let\mapsto=\longmapsto

\newcommand{\upperRomannumeral}[1]{\uppercase\expandafter{\romannumeral#1}}

\usepackage{tikz}
\usetikzlibrary{mindmap,backgrounds}
\usepackage{graphicx}
\usepackage[pagewise]{lineno}
\usepackage{stmaryrd}
\newtoggle{arxiv}
\togglefalse{arxiv}
\iftoggle{arxiv}{
	\usepackage{url}
	\bibliographystyle{plainurl}
	\typeout{In the arxiv mode}
	\usepackage{fancyhdr}
	\usepackage[colorlinks,citecolor=blue]{hyperref}
}
{
	\usepackage[backend=biber,style=alphabetic]{biblatex} 
	\usepackage{fancyhdr}
	\let\cite=\parencite
\addbibresource{../../master/master6.bib}
\addbibresource{hoshi-bib.bib}
\addbibresource{mochizuki-bib.bib}
\addbibresource{uchida-bib.bib}
\addbibresource{mochizuki-flowchart.bib}
}
\newtoggle{draft}
\toggletrue{draft}
\newtheorem{theorem}[equation]{Theorem}      
\newtheorem{lemma}[equation]{Lemma}          %
\newtheorem{corollary}[equation]{Corollary}  
\newtheorem{proposition}[equation]{Proposition}

\theoremstyle{definition}
\newtheorem{conj}[equation]{Conjecture}

\theoremstyle{definition}
\newtheorem{defn}[equation]{Definition}
\theoremstyle{remark}

\theoremstyle{definition}
\newtheorem{remark}[equation]{Remark}

\numberwithin{equation}{section}

\let\into=\hookrightarrow
\let\isom=\simeq

\newcommand{\A}{\mathcal{A}}

\newcommand{\bF}{{\bar{F}}}

\newcommand{\C}{{\mathbb C}}

\newcommand{\gal}{{\rm Gal}}

\newcommand{\N}{\mathcal{N}}

\newcommand{\Q}{{\mathbb Q}}

\renewcommand{\int}{\operatorname{int}}
\renewcommand{\O}{{\mathcal O}}
\renewcommand{\P}{{\mathbb P}}

\renewcommand{\bpro}{\begin{proposition}}
	\renewcommand{\epro}{\end{proposition}}
\renewcommand{\bcon}{\begin{conj}}
	\renewcommand{\econ}{\end{conj}}
\newcommand{\ilim}{\varprojlim}
\let\mathcal=\mathscr

\title{Untilts of fundamental groups: construction of labeled isomorphs of fundamental groups
}
\author{Kirti Joshi}

\begin{document}
\maketitle


\lhead{}

\epigraphwidth0.55\textwidth
\epigraph{Heard melodies are sweet, but those unheard \\ Are sweeter; therefore, ye soft pipes, play on;}{\iftoggle{arxiv}{John Keats \cite{keats}}{\citeauthor{keats}}}

\newcommand{\act}{\curvearrowright}
\newcommand{\lmp}{{\Pi\act\Ot}}
\newcommand{\lmpi}{{\lmp}_{\int}}
\newcommand{\lmpf}{\lmp_F}
\newcommand{\Om}{\O^{\times\mu}}
\newcommand{\Omf}{\O^{\times\mu}_{\bF}}
\renewcommand{\N}{\mathbb{N}}
\newcommand{\yoga}{Yoga}
\newcommand{\gl}[1]{{\rm GL}(#1)}
\newcommand{\bK}{\overline{K}}
\newcommand{\reptrip}{\rho:G_K\to\gl V}
\newcommand{\reptripp}[1]{\rho\circ\alpha:G_{\ifstrempty{#1}{K}{{#1}}}\to\gl V}
\newcommand{\benumlab}{\begin{enumerate}[label={{\bf(\arabic{*})}}]}
\newcommand{\ord}{\mathop{\rm ord}\nolimits}	
\newcommand{\kcs}{K^\circledast}
\newcommand{\lcs}{L^\circledast}
\renewcommand{\A}{\mathbb{A}}
\newcommand{\bfq}{\bar{\mathbb{F}}_q}
\newcommand{\tripod}{\P^1-\{0,1728,\infty\}}

\newcommand{\vseq}[2]{{#1}_1,\ldots,{#1}_{#2}}
\newcommand{\anab}[4]{\left({#1},\{#3 \}\right)\anabelmap\left({#2},\{#4 \}\right)}

\newcommand{\gln}{{\rm GL}_n}
\newcommand{\glo}[1]{{\rm GL}_1(#1)}
\newcommand{\glt}[1]{{\rm GL_2}(#1)}

\newcommand{\iut}{\cite{mochizuki-iut1, mochizuki-iut2, mochizuki-iut3,mochizuki-iut4}}
\newcommand{\topics}{\cite{mochizuki-topics1,mochizuki-topics2,mochizuki-topics3}}

\newcommand{\linv}{\mathfrak{L}}
\newcommand{\bedef}{\begin{defn}}
\newcommand{\eedef}{\end{defn}}
\renewcommand{\act}[1][]{\overset{#1}{\curvearrowright}}
\newcommand{\bfx}{\overline{F(X)}}
\newcommand{\anabelmap}{\leftrightsquigarrow}
\newcommand{\ban}[1][G]{\mathcal{B}({#1})}
\newcommand{\pit}{\Pi^{temp}}
 
 \newcommand{\bL}{\overline{L}}
 \newcommand{\bkm}{\bK_M}
 \newcommand{\vbk}{v_{\bK}}
 \newcommand{\vbkm}{v_{\bkm}}
\newcommand{\ocs}{\O^\circledast}
\newcommand{\ot}{\O^\triangleright}
\newcommand{\ocsk}{\ocs_K}
\newcommand{\otk}{\ot_K}
\newcommand{\ok}{\O_K}
\newcommand{\oko}{\O_K^1}
\newcommand{\oks}{\ok^*}
\newcommand{\Qpb}{\overline{\Q}_p}
\newcommand{\Qpbh}{\widehat{\overline{\Q}}_p}
\newcommand{\tr}{\triangleright}
\newcommand{\ocpt}{\O_{\C_p}^\tr}
\newcommand{\ocpf}{\O_{\C_p}^\flat}
\newcommand{\sG}{\mathscr{G}}
\newcommand{\sxfe}{\mathscr{X}_{F,E}}
\newcommand{\sxfep}{\mathscr{X}_{F,E'}}
\newcommand{\loglt}{\log_{\sG}}
\newcommand{\fc}{\mathfrak{t}}
\newcommand{\ku}{K_u}
\newcommand{\kup}{\ku'}
\newcommand{\kt}{\tilde{K}}
\newcommand{\sGpf}{\sG(\O_K)^{pf}}
\newcommand{\hgm}{\widehat{\mathbb{G}}_m}
\newcommand{\bE}{\overline{E}}

\tableofcontents
\togglefalse{draft}
\newcommand{\FF}{\cite{fargues-fontaine}}
\iftoggle{draft}{\pagewiselinenumbers}{\relax}

\section{Introduction}
I show that one can explicitly construct topologically/geometrically distinguishable data which provide isomorphic copies (i.e. \emph{isomorphs}) of the tempered fundamental group of a geometrically connected, smooth, quasi-projective variety over  $p$-adic fields. This is done via Theorem~\ref{thm:main} and Theorem~\ref{thm:main2}. Notably Theorem~\ref{thm:main2} also shows that the absolute Grothendieck conjecture fails for the class of Berkovich spaces (over algebraically closed perfectoid fields), arising as analytifications of geometrically connected, smooth, projective variety over  $p$-adic fields. 

The existence of distinctly labeled copies of the tempered fundamental groups is, as far as I understand, crucial to \iut, but produced in loc. cit. by entirely different means (for more on this labeling problem see Section~\ref{se:untilts-of-Pi}).  Let me also say at the onset that Mochizuki's Theory does not  consider passage to complete algebraically closed fields such as $\C_p$ and so my approach here is a significant point of departure from Mochizuki's Theory $\ldots$ and the methods of this paper do not use any results or ideas from Mochizuki's work. Nevertheless the results presented here establish unequivocally that  isomorphs of tempered (and \'etale) fundamental groups, of distinguishable provenance, exist and can be explicitly constructed.

The copies provided by Theorem~\ref{thm:main} and Theorem~\ref{thm:main2} arise from untilts of a fixed algebraically closed perfectoid field of characteristic $p>0$ and hence I call these copies \emph{untilts of fundamental groups}, or more precisely \emph{untilts of tempered fundamental groups}. 

An important consequence of these results is Corollary~\ref{cor:labeling-tempered}, which  provides a \emph{function from a suitable Fargues-Fontaine curve to the  isomorphism class of  the tempered fundamental group of a fixed variety (as above)} which provides a natural way of labeling the copies  obtained here by closed points of a suitable Fargues-Fontaine curve. 

In the last section of the paper I show that there is an entirely analogous theory of untilts of topological fundamental groups of connected Riemann surfaces.

This note began as a part of another note, \cite{joshi-log-link}, which I put into a limited circulation some time in July 2020, outlining my own  approach to some constructions of \iut. Peter Scholze immediately, but gently, pointed out that the section of \cite{joshi-log-link}, from which the present note is extracted,  needed some details. At that time I was  readying another note, \cite{joshi-gconj}, for wider circulation and addressing the issue noted by Scholze took longer and on the way I was able to  substantially strengthen and clarify my results (which appear here). So ultimately I decided that it would be best to publish the present note separately (while preparation of \cite{joshi-log-link} continued). My thanks are due to Peter Scholze,  and also to Yuichiro Hoshi, Emmanuel Lepage, and Jacob Stix, for  promptly providing comments, suggestions or corrections.

\newcommand{\ebh}{\widehat{\bE}}
\newcommand{\ebhx}[1][x]{\widehat{\bE^{#1}}}
\newcommand{\bdr}{B_{dR}}
\newcommand{\bdre}{{\bdr}_{,E}}
\newcommand{\bdrep}{B^+_{dR,E}}
\newcommand{\kbh}{\widehat{\bK}}

\section{The main theorem}

\newcommand{\bPi}{\overline{\Pi}}
\newcommand{\bPit}{\bPi^{\rm{\scriptscriptstyle temp}}}
\newcommand{\Pit}{\Pi^{\rm{\scriptscriptstyle temp}}}
\renewcommand{\pit}[1]{\Pi^{\scriptscriptstyle temp}_{#1}}
\newcommand{\pitk}[2]{\Pi^{\scriptscriptstyle temp}_{#1;#2}}
\newcommand{\pio}[1]{\pi_1({#1})}

\newcommand{\xan}{X^{an}}
\newcommand{\yan}{Y^{an}}

\emph{All valuations on base fields considered in this paper will be rank one valuations.} For the theory of tempered fundamental groups see \cite{andre-book,andre03} or \cite{lepage-thesis}. Berkovich spaces will be strictly analytic (and mostly will arise as analytifications of geometrically connected smooth quasi-projective varieties).

In what follows I will work with algebraically closed, perfectoid fields of characteristic zero. Such fields can also be characterized in many different ways.  For the convenience of the readers unfamiliar with perfectoid fields, the following lemma (immediate from \cite[Definition 3.1]{scholze12-perfectoid-ihes}), provides a translation of this condition into  more familiar hypothesis.

\blem\label{lem:perfectoid}
Let $K$ be a valued field and let $R\subset K$ be the valuation ring. The following conditions are equivalent:
\benumlab
\item $K$ is an algebraically closed field, complete with respect to a rank one non-archimedean valuation and with residue characteristic $p>0$.
\item $K$ is an algebraically closed, perfectoid field.
\eenum
\elem

\bp A perfectoid field has residue characteristic $p>0$ and is complete with respect to a rank one valuation. So  (2)$\implies$(1) is trivial. So it is enough to prove that (1)$\implies$(2).   I claim that  Frobenius $\phi:R/pR\to R/pR$ is surjective. let $\bar{x}\in R/pR$ and suppose $x\in R$ is an arbitrary lift of $\bar{x}$. Then as $K$ is algebraically closed, $x^{1/p}\mod{pR}$ provides a lift of $\bar{x}$. As $K$ is complete with respect to a rank one valuation and Frobenius is surjective on $R/pR$, so  $K$ is perfectoid by \cite[Definition 3.1]{scholze12-perfectoid-ihes} and by my hypothesis $K$ also algebraically closed. This proves (1)$\implies$(2).
\ep

For a perfectoid algebraically closed field $K$ as above, one has naturally associated field $K^\flat$, algebraically closed, perfectoid of characteristic $p>0$, called the \emph{tilt of $K$} and $K$ is called an untilt of $K^\flat$ (see \cite[Lemma 3.4]{scholze12-perfectoid-ihes}).

\newcommand{\cpt}{\C_p^\flat}
Fix an  algebraically closed  field, perfectoid $F$ of characteristic $p>0$ (see \cite{scholze12-perfectoid-ihes}). For example readers can simply assume, without any loss of generality, that $F= \C_p^\flat$ as this case is quite adequate for my purposes.

By an \emph{untilt} of $F$, I will mean a perfectoid field $K$, of characteristic zero, with $K^\flat$ isometric with $F$. Note that by \cite[Proposition 3.8]{scholze12-perfectoid-ihes}  $K$ is algebraically closed as its tilt $K^\flat=F$ is algebraically closed (by my hypothesis). If $F=\cpt$  then $K^\flat$ is isometric with $\cpt$. By the theory of \cite{fargues-fontaine} untilts $K$ of $F$ exist and are parametrized by Fargues-Fonatine curves.

Let $E$ be a $p$-adic field which is fixed for the present discussion. I will work with untilts $K$, of $F$, equipped with continuous embeddings $E\into K$ with the valuation of $K$ providing a valuation on $E$ which is equivalent to the natural $p$-adic valuation on $E$. By \cite{fargues-fontaine} for a given pair $(F,E)$,  such fields $K\hookleftarrow E$, exist and are parametrized by Fargues-Fontaine curves  (denoted here by $\sxfe$). Without further mention, all untilts $K$ will be assumed to be of this type (for our chosen $p$-adic field $E$).

Crucial point for this paper is that \emph{there exist untilts of $\cpt$ which are not topologically isomorphic.} This is the main result of \cite[Theorem~1.3]{kedlaya18}. Note that all characteristic zero untilts of $\cpt$ have the cardinality of $\C_p$ and are complete and algebraically closed fields and hence are abstractly isomorphic fields but may not be topologically isomorphic after \cite[Theorem 1.3]{kedlaya18}.

Now fix a geometrically connected, smooth quasi-projective variety $X/E$, with $E$ a $p$-adic field. Let $\xan/E$ be the strictly analytic space associated to $X/E$. Let \be\pit{X/E}=\pi_1^{temp}(\xan/E)\ee be the tempered fundamental group of the strictly $E$-analytic space associated to $X/E$ in the sense of \cite{andre03} or \cite{andre-book}. 

(\emph{Note that my notation $\pit{X/E}$ suppresses the passage to the   analytification $\xan/E$  for simplicity of notation. The theory of (tempered) fundamental groups also requires a choice of base point which will be suppressed from my notation.})

Let $E'/E$ be a finite extension of $E$ with a continuous embedding $E'\into K$ (as $K$ is algebraically closed, valued field containing $E$, such $E'$ exists.  Suppose $E'/E$ is also finite. One can consider $X_{E'}=X\times_{E}E'$ (similarly $X_K=X\times_EK$). Then one has an exact sequence by \cite[Prop. 2.1.8]{andre-book}
$$1\to \pit{X_{E'}/E'}\to \pit{X/E}\to \gal(E'/E)\to 1.$$

Let $\bE\subseteq K$ be the algebraic closure of $E$ contained in $K$.

By  varying $E'$ over all finite extensions of $E\into K$  one obtains (see \cite[Section 5.1]{andre03}) an exact sequence of topological groups:
$$1\to  \ilim_{E'/E}\pit{X_{E'}/E'}\to\pit{X/E} \to \gal(\bE/E)\to 1.$$

\newcommand{\Et}{\widetilde{E}}

\bthm\label{thm:main}
Let $F$ be an algebraically closed perfectoid field of characteristic $p>0$ (for example $F=\cpt$). Let $E$ be a $p$-adic field. Let $K,K_1,K_2$ be arbitrary untilts of $F$ with continuous embedding $E\into K$ (resp. into $K_1$ and $K_2$). Let $\bE$ (resp. $\bE_1,\bE_2$) be the algebraic closure of $E$ in $K$ (resp. in $K_1,K_2$). Let $X/E$ be a geometrically connected, smooth, quasi-projective variety over $E$. Then one has the following:
\benumlab
\item\label{thm:main-1} a continuous isomorphism $$\pit{X/K}\isom \ilim_{E'/E}\pit{X/E'},$$
where the inverse limit is over all finite extensions $E'$ of $E$ contained in $K$, and a
\item\label{thm:main-2} a short exact sequence of topological groups $$1\to \pit{X/K}\to \pit{X/E}\to G_E\to 1,$$ and
\item\label{thm:main-4} In particular for any two untilts $K_1,K_2$ of $F$, one has a continuous isomorphism $$\pit{X/K_1}\isom \pit{X/K_2}.$$
\eenum
\ethm
\bp 
The assertion \ref{thm:main-1} is the analogue of \cite[Prop. 5.1.1]{andre03} for an arbitrary untilt $K$ of $F$ containing $E$ (as above). Let me remind the reader that my hypothesis on  $K,K_1,K_2$ imply that $K,K_1,K_2$ are algebraically closed and complete with respect to a rank one valuation.  

Let me prove \ref{thm:main-1}, this will also lead to \ref{thm:main-2}. 
Since $K$ is  algebraically closed, it follows that $K$ contains an algebraic closure $\bE$ of $E$.  Let $\Et\subseteq K$ be the closure  (with respect to valuation topology of $K$) of $\bE$. 

It is clear that $\Et\supset \bE$ is complete and algebraically closed field and $\Et$ contains the algebraic closure $\bE\subset K$ of $E$ contained in $K$ as a dense subfield. In particular $\Et$ is the completion of $\bE$ with respect to the induced valuation. In other words $\Et$ is a copy of the completion of an algebraic closure of $E$ (usually denoted $\ebh$) equipped with an isometric embedding $\iota:\ebh\into K$ with $\iota(\ebh)=\Et$.
Hence  $K/\Et$ is an isometric extension of algebraically closed, complete valued fields (with rank one valuations). 

Now one can apply  the principle of invariance  of  fundamental groups under passage to extensions of algebraically closed fields. This principle is well-known for \'etale fundamental groups (see \cite[Expos\'e X, Corollaire 1.8]{sga1}). For  tempered fundamental groups this principle is proved in \cite[Proposition 2.3.2]{lepage-thesis}. Thus  applying \cite[Proposition 2.3.2]{lepage-thesis} to the extension $K/\Et$ one has an isomorphism of topological groups $$\pit{X_K}\isom \pit{X_{\Et}}.$$

On the other hand by \cite[Proposition 5.1.1]{andre03}, as $\Et$ is the completion of the algebraic closure of $\bE\subset K$ of $E$, one has an isomorphism
\be \pit{X/\Et}\isom \ilim_{E'/E}\pit{X_{E'}/E'}
\ee
and  an exact sequence of topological groups 
$$1\to \pit{X/K} \isom \ilim_{ E'/E}\pit{X_{E'}/E'}\to\pit{X/E} \to \gal(\bE/E)\to 1.$$
 This proves the assertions \ref{thm:main-1}, \ref{thm:main-2} as claimed.
 
Let me now prove \ref{thm:main-4}. The claimed isomorphism $\pit{X/K_1}\isom \pit{X/K_2}$ follows from the fact that both the groups can be identified with $\ilim_{E'/E}\pit{X/E'}$ where the inverse limit is over all finite extensions of $E'/E$ contained in $K_1$ (resp. $K_2$) and the fact that there is an equivalence between categories of finite extensions of $E$ contained in $K_1$ and the category of finite extensions of $E$ contained in $K_2$, since finite extensions of $E$ are given by adjoining roots of polynomials with coefficients in $E$ and this data is independent of the embedding of $E$ in $K_1$ or $K_2$ and moreover any abstract isomorphism of discretely valued fields is in fact an isometry--i.e given a finite extension of $E$,   $E'\into K_1$ contained in $K_1$, there is an isometry $E'\into K_2$ and vice versa.
\ep

\bthm\label{thm:main2}
Let $X/E$ be a geometrically connected, smooth projective variety. Let $K_1,K_2$ be two untilts of $\cpt$ which contain $E$. Suppose that $K_1,K_2$ are not topologically isomorphic. Then
\benumlab
\item\label{thm:main2-1} one has  an isomorphism of topological groups $$\pit{X/K_1}\isom \pit{X/K_2},$$
\item\label{thm:main2-2} but $\xan/K_1$ and $\xan/K_2$ are not isomorphic Berkovich spaces.
\item\label{thm:main2-3} In particular the absolute Grothendieck conjecture fails in the category of Berkovich spaces over perfectoid fields of characteristic zero.
\eenum
\ethm
\bp 
After Theorem~\ref{thm:main},
only \ref{thm:main2-2}  needs to be proved as \ref{thm:main2-2} $\implies$ \ref{thm:main2-3}. The hypothesis of Theorem~\ref{thm:main2} are non-vacuous--by\cite{kedlaya18}, fields $K_1,K_2$ exist. 

Assume that $X/E, K_1,K_2$ are as in my hypothesis and that $X$ is geometrically connected, smooth and projective over $E$. Suppose, if possible, that $\xan/K_1$ and $\xan/K_2$ are isomorphic as strictly analytic Berkovich spaces. Then one has a bounded isomorphism of Banach rings $$K_1\isom H^0(\xan/K_1,\O_{\xan/K_1}) \isom H^0(\xan/K_2,\O_{\xan/K_2})\isom K_2.$$ This isomorphism evidently extends  to a bounded  isomorphism of Banach fields $K_1\isom K_2$ and hence one sees that $K_1$ and $K_2$ are topologically isomorphic. Thus one has arrived at a contradiction.
\ep

\brem
As an aside let me remark that the proof of \cite{kedlaya18} (also see \cite[Th\'eor\`eme 2 and \S3 Remarque 2]{matignon84}) provides an uncountable collection of perfectoid fields $K_1,K_2$ with tilts isometric to $\cpt$ and such that  $K_1,K_2$ are not topologically isomorphic.
\erem

\section{Untilts of tempered fundamental groups}\label{se:untilts-of-Pi}
The results of the preceding section can be applied to the problem of producing \emph{labeled copies of the tempered fundamental groups.} A simple example of the labeling problem is the following: let $G$ be a topological group isomorphic to the absolute Galois group of some $p$-adic field. In this case one can ask if there are any distinguishable elements in the topological isomorphism class of $G$ with  the distinguishing features serving as labels.

For $G$ as above  the answer is simple: there is a distinguishable collection of copies of $G$, labeled by the $p$-adic fields $E$, i.e. $G_E\isom G$ as $p$-adic fields $E$, serving as labels  for copies of $G$ and the labels are distinguishable by their topological isomorphism class and so the label $E$ in $G_E$ is correspond to distinguishable geometric/topological datum of the $p$-adic field $E$. Moreover the main theorem of \cite{mochizuki-local-gro} also asserts that in fact the geometric label $K$ corresponds to an algebraic substructure of $G$ (``the upper numbering ramification filtration'') which provides the distinguishability. 

Now consider the labeling problem for the topological group $\Pi=\Pi_{X/E}$ for some hyperbolic curve $X$ over some $p$-adic field $E$. So one may again ask: is it possible to provide copies of $\Pi$ which are labeled by geometrically/topologically distinguishable labels?

As far as I understand this problem was considered, and solved (in many cases of interest) by Mochizuki in \iut.

\emph{Theorem~\ref{thm:main2} provides a different solution to this problem.} It shows that there exist isomorphic copies of tempered fundamental groups which arise from topologically distinct geometric data. The labels provided by these theorems arise from untilts of a fixed perfectoid, algebraically closed field of characteristic $p>0$ and so I call these copies of $\Pi$,  untilts of the fundamental group of $\Pi$.

By an \emph{untilt of the tempered fundamental group $\Pi=\pit{X/E}$} with respect to an untilt $K$ of $F$, I mean the tempered group $\pit{X/E}$ together with this short exact sequence
$$1\to \pit{X/K}\to \pit{X/E}\to G_E\to 1,$$
and I write $$\pitk{X/E}{K}$$ for this datum. Theorem~\ref{thm:main2} asserts that $\pit{X/K}$ is labelled by topologically and geometrically distinguished label $X/K$ especially as by Theorem~\ref{thm:main} and Theorem~\ref{thm:main2} one knows that if $K_1,K_2$ are two untilts of $F$ which are not topologically isomorphic (and hence non-isometric) then $K_1,K_2$ have inequivalent valuations so the spaces $\xan/K_1$ and $\xan/K_2$  are not be isomorphic rigid analytic spaces. So while $\pitk{X/E}{K_1}\isom \pitk{X/E}{K_2}$ are isomorphic as tempered groups, these arise from possibly distinct geometric spaces. So the terminology of untilting makes sense. 

Note that the labeling also provides a algebraic substructure of $\Pi$ namely the normal subgroup $\pit{X/K}\subset \pit{X/E}$.

In particular one has the following corollary:

\bcor\label{cor:labeling-tempered}
Let $X/E$ be a geometrically connected, smooth quasi-projective variety over a $p$-adic field $E$.
Then the natural function $$K\mapsto \pitk{X/E}{K}$$ from the set of inequivalent untilts of $\C_p^\flat$ to the topological isomorphism class of $\pit{X/E}$
provides a distinguished collection of distinctly labeled copies  $$\left\{\pitk{X/E}{K_x}: x \in \sxfe \text{ a closed point  with residue field } K_x \right\}$$ of the tempered fundamental group  $\pit{X/E}$.
\ecor

The above consideration can be applied to \'etale fundamental groups of geometrically connected, smooth quasi-projective varieties as follows. Let $X/E$ be a geometrically connected, smooth, quasi-projective variety over a $p$-adic field $E$. Then one has a natural homomorphism (\cite[Proposition 4.4.1]{andre03}, \cite[Section 2.1.4]{andre-book}):
$$\pit{X/E}\to\pi_1(X/E),$$ 
which is injective if $\dim(X)=1$, and in any dimension
its image is dense and moreover  $\pi_1(X/E)$ is the profinite completion  
$$\widehat{\pit{X/E}}=\pi_1(X/E).$$ 

Let $K$ be an untilt of $F$. I define
$$\pi_1(X/E)_{K}=\widehat{\pitk{X/E}{K}},$$
and call $\pi_1(X/E)_{K}$ the \emph{untilt of the \'etale fundamental group} $\pi_1(X/E)$ corresponding to the untilt $K$ (of $F\isom K^\flat$). 
Thus one has the notion of untilts of $\pi_1(X/E)$.

\bcor
Let $X/E$ be a geometrically connected, smooth quasi-projective variety over a $p$-adic field $E$.
Then the natural function $$K\mapsto \pi_1(X/E)_{K}$$ from the set of inequivalent untilts of $F$ to the topological isomorphism class of $\pi_1(X/E)$
provides a distinguished collection of distinctly labeled copies  $$\left\{\pi_1(X/E)_{K_x}: x \in \sxfe \text{ a closed point  with residue field } K_x \right\}$$ of the \'etale fundamental group  $\pi_1(X/E)$.
\ecor

I have used perfectoid algebraically closed fields as a set of distinguishing labels for the copies of fundamental groups produced here. There is in fact a bigger indexing set:

\newcommand{\sK}{\mathcal{K}}

\bcor
Let $E$ be a $p$-adic field, $X/E$ a geometrically connected, smooth, quasi-projective variety over $E$. Consider the set of topological isomorphism classes of algebraically closed, complete valued fields $K\supset E$ (isometric inclusions):  $$\sK_E=\left\{K: E\subset K, K=\kbh\right\}.$$

Then there is a natural function $K\mapsto \pi_1(X/E)_K$ from $\sK_E$ to the topological isomorphism class of the profinite group $\pi_1(X/E)$ given by considering the tempered fundamental group associated to the datum $(X,E\into K)$.
\ecor

\brem\label{re:archimedean-case}
There is a further aspect of this result which should be pointed out. One should view elements $K\in\sK_E$ as providing a topological variation of ambient (additive) structure   $K\supset \bE$ while keeping internal field structure of $\bE$ unchanged. \emph{Such variations exist because, unlike the number field case, $p$-adic fields, even complete algebraically closed fields such as $\C_p$, are quite far from being topologically rigid.} 
This is in complete contrast with the archimedean case, where by the well-known theorem of Ostrowski \cite[Chap. 6, \S 6, Th\'eor\`eme 2]{bourbaki-alg-com}, one knows that the only algebraically closed field complete with respect to an archimedean valuation is isometric to $\C$. To put Ostrowski's Theorem differently: \emph{Any two  algebraically closed, \emph{archimedean perfectoid fields} (i.e. fields which are algebraically closed and complete with respect to an archimedean valuation) are isometric (and also isometric with $\C$) and hence such fields are topologically rigid.}
\erem

\section{Untilts of fundamental groups of Riemann surfaces}
Let me point out that there is a complex analytic  analogue of the theory of untilting which is outlined above. Let $\Pi=\pi_1^{top}(X)$ be the topological fundamental group of a connected Riemann surface $X$, which one assumes to be hyperbolic to avoid trivialities. Then consider all connected Riemann surfaces $X'$ whose topological fundamental group $\pi_1^{top}(X')\isom \Pi$. One can also fix the genus $g$  and number of punctures $n$ of $X$ and consider $X'$ of genus $g$ with $n$ punctures. 

The assignment $X'\mapsto \pi_1^{top}(X')\isom \Pi$ provides a function from the isomorphism classes of connected, hyperbolic Riemann surfaces of genus $g$ and with $n$ punctures to the isomorphism class of the group $\Pi$. Then $\pi_1^{top}(X')$ is an untilt of $\Pi$ with the complex structure of $X'$ serving as a geometrically distinguishable feature of this copy of $\Pi$. 

Now assume that $K$ is a number field and $X/K$ is a hyperbolic, geometrically connected smooth quasi-projective curve. In \cite{tamagawa97-gconj}, \cite{mochizuki96-gconj} it has been shown that the genus $g$ of $X$ and the number of punctures on $X$ is amphoric i.e. determined by the isomorphism class of the topological group $\pi_1^{\acute{e}t}(X/K)$. So one can fix $g,n$. 

An \emph{untilt of $\pi_1^{\acute{e}t}(X)$ at $\infty$} (here $\infty$ is short for ``\emph{at  archimedean primes}'') is a pair consisting of an embedding $K\into \C$ and a  Riemann surface $X'$, of genus $g$ and with  $n$ punctures, such that $\widehat{\pi_1^{top}(X')}\isom \pi_1(X/\bK)$, where $\hat{-}$ denotes the profinite completion. An untilt of $\Pi=\pi_1(X/K)$ at $\infty$ will be labeled $\Pi_{K\into\C,X'}$. Two untilts of $\pi_1(X/K)$ at $\infty$ are equivalent if the the two embeddings of $K\into \C$ are equivalent (in the obvious sense) and the two corresponding Riemann surfaces are isomorphic. 

\newcommand{\cM}{\mathcal{M}}
Thus one has the following tautology:
\bpro
Fix a profinite group $\Pi\isom \pi_1(X/K)$ with $X/K$ a geometrically connected, smooth, hyperbolic curve over a number field $K$ with no real embeddings. Then the equivalence classes of untilts of $\Pi$ at $\infty$ are in bijection with $$\widetilde{Hom(K,\C)}\times \cM_{g,n},$$ where $\widetilde{Hom(K,\C)}$ is the set of equivalence classes of embeddings of $K \into \C$ and $\cM_{g,n}$ is the moduli stack of Riemann surfaces of genus $g$ with $n$ punctures.
\epro

\brem
Owing to the topological rigidity of algebraically closed fields complete with respect to an archimedean absolute value, forced by Ostrowski's Theorem (see Remark \ref{re:archimedean-case}), one could say that  untilts of topological fundamental groups at $\infty$ (i.e. at archimedean primes) can  arise only from the existence of geometric variations of the underlying objects. 
\erem

\brem 
Readers familiar with the classical Szpiro inequality (for surfaces fibered over curves) and its several  different proofs, may notice that the above proposition provides a unified way of viewing these proofs  as taking place  over the ``space of untilts.'' More precisely the ``space of untilts'' provides the geometric Kodaira-Spencer classes which underly these proofs. 
\erem

\iftoggle{arxiv}{
	\bibliography{hoshi-bib,mochizuki-bib,uchida-bib,mochizuki-flowchart,../../master/master6.bib}
}
{
	\printbibliography
}

\end{document}